\documentclass[11pt,final]{article}%
\usepackage{stix2}
\usepackage{amsfonts}
\usepackage{amsmath}
\usepackage{amsxtra}
\usepackage[amsmath,hyperref,thmmarks]{ntheorem}
\usepackage{graphicx}
\usepackage[colorlinks=true,bookmarksnumbered=true,pdfpagemode=None,final]%
{hyperref}
\usepackage{amssymb}
\usepackage{enotez}%
\providecommand{\U}[1]{\protect\rule{.1in}{.1in}}

\newcommand{\df}{\smash{\lower.12em\hbox{\textup{\tiny def}}}}

\usepackage[overload]{textcase}

\usepackage{xcolor}
\definecolor{bblue}{rgb}{0.0, 0.0, 0.6}

\usepackage{array}

\usepackage{microtype}

\usepackage{tikz}
\usepackage{tikz-cd}
\usetikzlibrary{matrix,arrows,positioning,decorations.pathmorphing}
\usetikzlibrary{decorations.markings,shapes.geometric}
\tikzset{commutative diagrams/column sep/Huge/.initial=24ex}
\tikzcdset{arrow style=math font}

\usepackage{wrapfig}

\usepackage{enumitem}
\setitemize[1]{label=$\diamond$\hspace{0.07in}}
\setlist{topsep=0.2em,itemsep=0.2em,parsep=0.2em}

\setdescription{font=\normalfont}

\usepackage{url}
\usepackage{verbatim}

\usepackage[notcite,notref]{showkeys}

\usepackage{mathtools}

\newcommand{\eb}[1]{{\itshape\bfseries#1}}
\renewcommand{\emph}{\eb}

\usepackage{titletoc,titlesec}
\contentsmargin{2.0em}
\dottedcontents{section}[2.3em]{}{1.8em}{1pc}
\usepackage[nottoc]{tocbibind}

\titleformat*{\section}{\LARGE\bfseries}
\titleformat*{\subsection}{\Large\itshape}
\titleformat*{\subsubsection}{\scshape}
\titleformat*{\paragraph}{\itshape}

\usepackage{natbib}
\let\cite\citealt

\hypersetup{linkcolor=bblue,anchorcolor=bblue,citecolor=bblue,filecolor=bblue,urlcolor=bblue}
\hypersetup{pdftitle={A Primer of Commutative Algebra},pdfauthor={J.S. Milne}}
\hypersetup{pdfsubject={commutative algebra},pdfkeywords={ideals}}

\usepackage[papersize={6.5in,10.0in},margin=0.5in,pdftex]{geometry}

\DeclareMathOperator{\vect}{mod}

\newcommand{\babstract}{\begin{abstract}}\newcommand\eabstract{\end{abstract}}
\newcommand{\bcomment}{}
\newcommand{\bfootnotesize}{\begin{footnotesize}}\newcommand\efootnotesize{\end{footnotesize}}
\newcommand{\bquote}{\begin{quote}}\newcommand\equote{\end{quote}}
\newcommand{\bsmall}{\begin{small}}\newcommand\esmall{\end{small}}
\newcommand{\btable}{\begin{table}}\newcommand{\etable}{\end{table}}

\newcommand{\edocument}{\end{document}}

\newcommand{\tstyle}{\textstyle}

\def\1{{1\mkern-7mu1}}

\DeclareMathOperator{\Aut}{Aut}

\DeclareMathOperator{\End}{End}

\DeclareMathOperator{\GL}{GL}

\DeclareMathOperator{\Hom}{Hom}
\DeclareMathOperator{\id}{id}

\DeclareMathOperator{\im}{Im}  

\DeclareMathOperator{\Ker}{Ker}
\DeclareMathOperator{\Lie}{Lie}

\DeclareMathOperator{\rank}{rank}

\DeclareMathOperator{\sign}{sign}
\DeclareMathOperator{\SL}{SL}

\DeclareMathOperator{\Spec}{Spec}

\DeclareMathOperator{\Sym}{Sym}

\newcommand{\Rep}{\mathsf{Rep}}

\theoremnumbering{arabic}
\theoremheaderfont{\scshape}
\RequirePackage{latexsym}
\theorembodyfont{\slshape}
\theoremseparator{.}
\newtheorem{X}{X}[section]

\newtheorem{lemma}{Lemma}
\newtheorem{proposition}{Proposition}
\newtheorem{theorem}{Theorem}
\newtheorem{corollary}{Corollary}
\theorembodyfont{\upshape}

\newtheorem{plain}[X]{}

\newtheorem{question}{Question}
\newtheorem{remark}{Remark}

\theorembodyfont{\small}

\theorembodyfont{\normalsize}
\theoremstyle{nonumberplain}
\theoremsymbol{\ensuremath{\color{lightgray}\blacksquare}}
\newtheorem{proof}{Proof}
\qedsymbol{\ensuremath{\color{lightgray}\blacksquare}}
\begin{document}

\title{Algebraic groups as automorphism groups of algebras}
\author{J.S. Milne}
\date{December 23, 2022}
\maketitle

\begin{abstract}
We show that every algebraic group scheme over a field with at least 8
elements can be realized as the group of automorphisms of a nonassociative
algebra. This is only a modest improvement of the theorem of Gordeev and Popov
(2003), but it allows us to give a new characterization of algebraic Lie
algebras and to simplify the standard descriptions of Mumford--Tate domains
and Shimura varieties as moduli spaces. Once the original argument of Gordeev
and Popov has been rewritten in the language of schemes, we find that it also
applies to algebraic groups over Dedekind domains.

\end{abstract}
\tableofcontents

\section*{Introduction}

\addcontentsline{toc}{section}{Introduction}

Let $k$ be a field. We use the following conventions: an \textit{algebra }%
$A$\textit{ over }$k$ is a $k$-vector space $V$ equipped with a $k$-linear map
$t\colon V\otimes_{k}V\rightarrow V$ (no conditions);\footnote{This is
Bourbaki's definition. Note that we do not require an algebra to have a
two-sided identity.} a \textit{commutative }$k$\textit{-algebra} is a
commutative associative $k$-algebra with an identity element; an
\textit{algebraic group over }$k$\ is an affine group scheme of finite type
over $k$. When $V$ is a vector space over $k$ and $R$ is a commutative
$k$-algebra, $V_{R}$ denotes the $R$-module $R\otimes_{k}V$.

Let $A=(V,t)$ be a finite-dimensional algebra over $k$. The functor
$R\rightsquigarrow\Aut_{R}(A_{R})$ of commutative $k$-algebras is represented
by an algebraic subgroup of $\GL_{V}$, which we denote by $\underline{\Aut}%
(A)$. In general, $\underline{\Aut}(A)$ need not be smooth.

In the remainder of the introduction, $k$ is a field with at least $8$ elements.

\begin{theorem}
\label{T1}Every algebraic group over $k$ is isomorphic to $\underline{\Aut}%
(A)$\, for some finite-dimensional algebra $A$ over $k$.
\end{theorem}

This statement is almost the same as that of Theorem 1 in Gordeev and Popov
2003. However, there \textquotedblleft algebraic group\textquotedblright\ is
meant in the sense of Borel 1991, not schemes. In the language of schemes,
they prove that, for each smooth algebraic group $G$ over $k$, there exists a
finite-dimensional algebra $A$ over $k$ such that $G(K)=\Aut(K\otimes A)$ for
all fields $K$ containing $k$ (ibid., Corollary 1). By contrast, we prove
that, for each algebraic group $G$ (not necessarily smooth) over $k$, there
exists a finite-dimensional algebra $A$ over $k$ such that $G(R)=\Aut(R\otimes
A)$ for all commutative $k$-algebras $R$.

If $G=\underline{\Aut}(A)$, then, in particular, $G(R)=\Aut(R\otimes A)$ for
$R$ the ring of dual numbers over $k$. From this it follows that $\Lie(G)$ is
the Lie algebra $\mathrm{Der}(A)$ of derivations of $A$. We now have the
following simple criterion for a Lie algebra to be algebraic.

\begin{corollary}
\label{C11}Let $\mathfrak{g}{}$ be a finite-dimensional Lie algebra over $k$.
Then $\mathfrak{g}{}=\Lie(G)$ for some algebraic group $G$ over $k$ if and
only if $\mathfrak{g}{}=\mathrm{Der}(A)$ for some finite-dimensional algebra
$A$ over $k$.
\end{corollary}

\begin{proof}
If $\mathfrak{g}{}=\mathrm{Der}(A)$, we define $G$ to be $\underline{\Aut}%
(A)$, and then $\Lie(G)=\mathrm{Der}(A)=\mathfrak{g}{}$. Conversely, if
$\mathfrak{g}{}=\Lie(G)$, we use Theorem \ref{T1} to set $G=\underline{\Aut}%
(A)$, and then $\mathrm{Der}(A)=\Lie(G)=\mathfrak{g}$.
\end{proof}

\begin{remark}
When $k$ has characteristic $p\neq0$, both $\Lie(G)$ and $\mathrm{Der}(A)$
have natural $p$-Lie algebra structures, and Corollary \ref{C11} holds with
\textquotedblleft Lie algebra\textquotedblright\ replaced by \textquotedblleft%
$p$-Lie algebra\textquotedblright.
\end{remark}

Theorem \ref{T1} extends to neutral tannakian categories. An \textit{algebra}
in a tensor category is an object $X$ equipped with an algebra structure,
i.e., a morphism $t\colon X\otimes X\rightarrow X$.

\begin{corollary}
\label{C12}Let $\mathsf{C}$ be a neutral algebraic\footnote{A tannakian
category over a field is said to be algebraic if it corresponds to an
algebraic gerbe. This amounts to saying that the affine group scheme attached
to a fibre functor over some extension field of the base field is algebraic,
i.e., of finite type. See Saavedra 1972, III, 3.3.1.} tannakian category over
$k$. There exists an algebra $(X,t)$ in $\mathsf{C}$ such that, for every
fibre functor $\omega$ with values in a field $k^{\prime}\supset k$,%
\[
\underline{\Aut}^{\otimes}(\omega)=\underline{\Aut}(\omega(X),\omega
(t))\text{.}%
\]

\end{corollary}

\begin{proof}
As $\mathsf{C}$ is neutral, there exists a $k$-valued fibre functor
$\omega_{0}$, and $\omega_{0}$ defines an equivalence of tensor categories
$\mathsf{C}\rightarrow\Rep(G)$, where $G=\underline{\Aut}^{\otimes}(\omega
_{0})$. According to Theorem 1, $G=\underline{\Aut}(A)$ for some algebra
$A=(V,t)$ in $\Rep(G)$. There exists an algebra $(X,t)$ in $\mathsf{C}$ such
that $\omega_{0}(X,t)$ is isomorphic to $(V,t)$. For any $k^{\prime}$-valued
fibre functor $\omega$,
\[
\underline{\Aut}^{\otimes}(\omega)\subset\underline{\Aut}(\omega
(X),\omega(t)),
\]
but $\omega$ becomes isomorphic to $\omega_{0}$ over some field containing
$k^{\prime}$, and so the inclusion is an equality.
\end{proof}

\begin{question}
\label{Q1}Does Corollary \ref{C12} hold for nonneutral tannakian categories?
\end{question}

Let $G$ be an algebraic group over $k$. A standard result says that there
exists a finite-dimensional $k$-vector space $V$ and a family of tensors for
$V$ such that $G$ is isomorphic to the subgroup of $\GL_{V}$ fixing the
tensors. Theorem \ref{T1} gives a more precise statement.

\begin{corollary}
\label{C13}Let $G$ be an algebraic group over $k$. There exists a
finite-dimensional $k$-vector space $V$ and a $t\in V\otimes V^{\vee}\otimes
V^{\vee}$ such that $G$ is isomorphic to the subgroup of $\GL_{V}$ fixing $t$.
Here $V^{\vee}$ is the linear dual of $V$.
\end{corollary}

\begin{proof}
Let $V$ be a finite-dimensional $k$-vector space, and let $t^{\prime}\colon
V\otimes V\rightarrow V$ be the linear map corresponding to $t\in V\otimes
V^{\vee}\otimes V^{\vee}$. Let $R$ be a commutative $k$-algebra and $\alpha$
an $R$-linear automorphism of $V_{R}$. Then $\alpha(t)=t$ if and only if
$\alpha$ is an algebra automorphism of $(V,t^{\prime})$. Thus, the corollary
is a restatement of Theorem \ref{T1}.
\end{proof}

Once the proof of Theorem 1 of Gordeev and Popov has been rewritten in the
language of schemes, one sees that it in fact applies over more general bases.
In particular, we prove the following statement.

\begin{theorem}
\label{T3}Let $G$ be a flat algebraic group over a Dedekind domain $R$. If $R$
has enough units, then there exists an algebra $A$ over $R$, flat and finitely
generated as an $R$-module, such that $G$ is isomorphic to $\underline{\Aut}%
(A)$ (i.e., $G$ represents the functor of commutative $R$-algebras $R^{\prime
}\rightsquigarrow\Aut(A_{R^{\prime}})$).
\end{theorem}

See Theorem \ref{a36} for a precise statement. In the course of proving
Theorem \ref{T3}, we obtain the following result (Corollary \ref{a41} to
Proposition \ref{a40}).\footnote{The second statement should be considered
folklore. Gabber has proved a similar result over noetherian regular base
schemes of dimension $\leq2$ --- see the revised 2011 version of SGA 3,
Expos\'{e} VI$_{B}$, Prop.\, 13.2.}

\begin{theorem}
\label{T2}Let $G$ be a flat algebraic group over a Dedekind domain $R$. There
exists a finite flat $R$-submodule $V$ of $\mathcal{O}{}(G)$, stable under
$G$, such that the homomorphism $G\rightarrow\GL_{V}$ is a closed immersion.
If $R$ is principal and the generic fibre of $G$ over $R$ is linearly
reductive, then $G$ is the subgroup of $\GL_{V}$ fixing a finite collection of
tensors in spaces $V^{\otimes m}\otimes(V^{\vee})^{\otimes n}$.
\end{theorem}

As polarizable rational Hodge structures form a tannakian category${}$, it is
possible to equip such a Hodge structure with an algebra structure. Theorem
\ref{T1} allows us to realize Mumford--Tate domains as a moduli spaces for
polarized rational Hodge structures with an algebra structure. This is simpler
than the usual description in terms of polarized rational Hodge structures
equipped with some family of Hodge tensors. 

Similarly, Theorem 1 and its corollaries allow us to realize Shimura varieties
of abelian type with rational weight as moduli schemes for abelian motives
with an algebra structure. This is simpler than the description in Theorems
3.13 and 3.31 of  Milne 1994. As this depends on Deligne's theory of absolute
Hodge classes on abelian varieties, at present it applies only in
characteristic zero. However, once Deligne's theory has been extended to mixed
characteristic (cf.~Milne 2009), Theorem 2 will allow us to obtain a new
moduli description of Shimura varieties in mixed characteristic. This should
allow a significant simplification of the theory. It was this that sparked the
author's interest in the topic. We do not explain these applications here as we
plan to return to them in a future work. For a brief explanation, see the last
two sections of arXiv:2012.05708v1.

In Section 1 of the article we explain, following Gordeev and Popov, the
construction of some algebras, and in Section 2 we prove our main theorems.

\subsection*{Notation and conventions.}

Throughout, $R$ is a commutative ring with $1$. By a \textit{finite flat}
$R$-module, we mean an $R$-module that is finitely presented and flat.
Unadorned tensor products are over $R$. We say that an $R$-module $S$ is a
\textit{direct summand} of an $R$-module $W$ if it is a submodule of $W$ and
admits a complement, i.e., $W=S\oplus W^{\prime}$ for some $W^{\prime}$.

An \textit{algebra} $A$ over $R$ is an $R$-module $V=\mathrm{mod}(A)$ together
with an $R$-linear map $t\colon V\otimes_{k}V\rightarrow V$. We say that $A$
is finitely presented, flat, \ldots if the $R$-module $\mathrm{mod}(A)$ is
finitely presented, flat, \ldots. For an element $a$ of an algebra, $r_{a}$
denotes right multiplication by $a$. We let $\langle S\rangle$ denote the
linear span of a subset $S$ of a module.

For a finite flat $R$-module $V$, we let $\mathrm{T}(V)$ denote the tensor
algebra of $V$,
\[
\mathrm{T}(V)=\bigoplus\nolimits_{i\geqslant0}V^{\otimes i},
\]
and we let $\mathrm{T}(V)_{+}$ denote the following ideal in $T(V)$,%
\begin{equation}
\mathrm{T}(V)_{+}=\bigoplus\nolimits_{i\geqslant1}V^{\otimes i}\text{.}
\nonumber\label{1.1}%
\end{equation}
Both are algebras over $R$ equipped with a natural action of the algebraic
group $\GL_{V}$:
\[
g\cdot t_{i}=g^{\otimes i}(t_{i}),\quad g\!\in\!\GL(V_{R}),\quad t_{i}\in
V_{R}^{\otimes i},\quad R\text{ a commutative }k\text{-algebra.}%
\]

By an algebraic group over $R$, we mean an affine group scheme of finite
presentation over $R$. An embedding of algebraic groups is a morphism that is
both a homomorphism and a closed immersion. For an $R$-module $V$ with an
action of an algebraic group $G$, we let $V_{0}$ denote $V$ equipped with the
trivial action of $G$.

\section{Some special algebras}

This section is adapted from Gordeev and Popov 2003 and Perepechko 2009.

\subsection{The algebra $A(V,S)$}

\begin{proposition}
\label{S21} Let $V$ be a nonzero finite flat $R$-module. Let $S$ be a finite
flat $R$-submodule of $V^{\otimes r}$, some $r>1$, such that $V^{\otimes r}/S$
is flat. Then there exists a finite flat graded algebra\footnote{Here $A^{2}$
is the part of degree $2$ of the graded algebra $A$.} $A=V\oplus A^{2}\oplus\cdots$
over $R$ such that
\[
(\GL_{V})_{S}=\underline{\Aut}(A,V)\qquad\text{(automorphisms of $A$
stabilizing $V$).}%
\]

\end{proposition}

Here $(\GL_{V})_{S}$ represents the functor $R^{\prime}\rightsquigarrow
\{g\in\GL(V_{R^{\prime}})\mid g^{\otimes r}(S_{R^{\prime}})=S_{R^{\prime}}\}.$

\begin{proof}
Let
\[
I(S)=S\oplus\Big(\bigoplus\nolimits_{i>r}V^{\otimes i}\Big).
\]
It is an ideal in the algebra $\mathrm{T}(V)_{+}$, and we define
\[
A(V,S)=\mathrm{T}(V)_{+}/I(S).
\]
This is a finite flat algebra over $R$ with%
\[
\vect(A(V,S))=\Big(\bigoplus\nolimits_{i=1}^{r-1}V^{\otimes i}\Big)\oplus
(V^{\otimes r}/S)
\]
as a graded $R$-module. Let $R^{\prime}$ be a commutative $R$-algebra. When we
replace $V$ and $S$ with $V_{R^{\prime}}$ and $S_{R^{\prime}}$ in the above
definition, we obtain an algebra $A(V_{R^{\prime}},S_{R^{\prime}})$ over
$R^{\prime}$, and%
\[
A(V_{R^{\prime}},S_{R^{\prime}})\simeq R^{\prime}\otimes_{R}A(V,S).
\]

The ideal $I(S)$ is stable under the natural action of $(\GL_{V})_{S}$ on
$\mathrm{T}(V)_{+}$, and so $(\GL_{V})_{S}$ acts on $A(V,S)$ by algebra
automorphisms. The quotient map $\pi\colon\mathrm{T}(V)_{+}\rightarrow A(V,S)$
is $(\GL_{V})_{S}$-equivariant. The condition $r>1$ ensures that
$V=V^{\otimes1}$ is a submodule of $A(V,S)$. Hence $(\GL_{V})_{S}$ acts
faithfully on $A(V,S)$, and it stabilizes $V$. It remains to show that the
algebraic group $(\GL_{V})_{S}$ represents the functor%
\[
R^{\prime}\rightsquigarrow\{\sigma\in\Aut(A(V,S)_{R^{\prime}})\mid
\sigma(V_{R^{\prime}})=V_{R^{\prime}}\}.
\]

Let $R^{\prime}$ be a commutative $R$-algebra. We have seen that%
\[
(\GL_{V})_{S}(R^{\prime})\subset\{\sigma\in\Aut(A(V_{R^{\prime}},S_{R^{\prime
}}))\mid\sigma(V_{R^{\prime}})=V_{R^{\prime}}\}
\]
and it remains to prove equality. Let $\sigma$ be an element
$\Aut(A(V_{R^{\prime}},S_{R^{\prime}}))$ such that $\sigma(V_{R^{\prime}%
})=V_{R^{\prime}}$. Put $g=\sigma|_{V_{R^{\prime}}}$, and let $g^{\bullet}$
denote the canonical extension of $g$ to an automorphism of $\mathrm{T}%
(V_{R^{\prime}})_{+}$. Then $g^{\bullet}|_{V_{R^{\prime}}}=g=\sigma
|_{V_{R^{\prime}}}$, and the diagram
\[
\begin{tikzcd}
\mathrm{T}(V_R)_+\arrow{r}{g^{\bullet}}\arrow{d}{\pi_R}
&\mathrm{T}(V_R)_+\arrow{d}{\pi_R}\\
A(V_R,S_R)\arrow{r}{\sigma}&A(V_R,S_R)
\end{tikzcd}
\]
commutes because it does on $V_{R^{\prime}}$, which generates the algebra
$\mathrm{T}(V_{R^{\prime}})_{+}$. The commutativity of the diagram implies
that $g^{\bullet}(\Ker(\pi_{R^{\prime}}))=\Ker(\pi_{R^{\prime}})$. As
$\Ker(\pi_{R^{\prime}})=I(S_{R^{\prime}})$, it follows that $g$ is an element
of $\GL(V_{R^{\prime}})$ such that $g^{\otimes r}(S_{R^{\prime}}%
)=S_{R^{\prime}}$. The diagram shows that its image in $\Aut(A(V_{R^{\prime}%
},S_{R^{\prime}}))$ is $\sigma$.
\end{proof}

\subsection{Two lemmas}

\label{S22}

\begin{lemma}
\label{L0}Let $V$ be a finite flat $R$-module and $\phi$ an endomorphism of
$V$. Suppose that $V$ decomposes into a direct sum of eigenspaces%
\[
V=V_{1}\oplus\cdots\oplus V_{n}%
\]
for $\phi$ with eigenvalues $\alpha_{1},\ldots,\alpha_{n}\in R$ that are
distinct modulo every maximal ideal of $R$. For any commutative $R$-algebra
$R^{\prime}$, $V_{R^{\prime}}=V_{1R^{\prime}}\oplus\cdots\oplus V_{nR^{\prime
}}$ with
\begin{equation}
V_{iR^{\prime}}=\{x\in V_{R^{\prime}}\mid\phi_{R^{\prime}}(x)=\alpha_{i}x\}.
\label{q2}%
\end{equation}

\end{lemma}

\begin{proof}
Certainly, $V_{R^{\prime}}\overset{\df}{=}R^{\prime}\otimes_{R}V$ is the
direct sum of the $R^{\prime}$-modules $V_{iR^{\prime}}\overset{\df}{=}%
R^{\prime}\otimes_{R}V_{i}$ and $V_{iR^{\prime}}$ is contained in the
right-hand side of (\ref{q2}). For the opposite inclusion, let $x\in
V_{R^{\prime}}$ be such that $\phi_{R^{\prime}}(x)=\alpha_{i}x$, and write
$x=x_{1}+\cdots+x_{n}$ with $x_{j}\in V_{jR^{\prime}}$. Then%
\[
\phi_{R^{\prime}}(x)=\alpha_{1}x_{1}+\cdots+\alpha_{n}x_{n}%
\]
and so
\[
\textstyle0=\phi_{R^{\prime}}(x)-\alpha_{i}x=\sum_{j}(a_{j}-\alpha_{i}%
)x_{j}\text{.}%
\]
It follows that $(a_{j}-\alpha_{i})x_{j}=0$ for all $j\neq i$. As $(\alpha
_{j}-\alpha_{i})\in R^{\times}\subset R^{\prime\times}$, this implies that
$x_{j}=0$ for all $j\neq i$.
\end{proof}

\begin{lemma}
\label{L1} Let $A$ be an algebra over $R$ with a left identity element $e\in
A$. Suppose that $\vect(A)$ decomposes into a direct sum of eigenspaces
\[
\vect(A)=Re\oplus A_{1}\oplus{\cdots}\oplus A_{r}%
\]
for $r_{e}$ with eigenvalues $1,\alpha_{1},\ldots,\alpha_{r}\in R$ such that
$0,1,$ $\alpha_{1}$, $\ldots,\alpha_{r}$ are distinct modulo every maximal
ideal of $R$. For any commutative $R$-algebran $R^{\prime}$,

\begin{enumerate}
\item $e$ is the unique left identity element in $A_{R^{\prime}}$;

\item if $\sigma\in\Aut(A_{R^{\prime}})$\textrm{,} then $\sigma(e)=e$ and
$\sigma(A_{iR^{\prime}})=A_{iR^{\prime}}$ for all $i$.
\end{enumerate}
\end{lemma}

\begin{proof}
According to Lemma \ref{L0},
\[
\vect(A_{R^{\prime}})=R^{\prime}e\oplus A_{1R^{\prime}}\oplus{\cdots}\oplus
A_{rR^{\prime}}%
\]
with $A_{iR^{\prime}}=\{x\in A_{R^{\prime}}\mid xe=\alpha_{i}x\}$.

(a) Let $e^{\prime}$ be a left identity element of $A_{R^{\prime}}$, and write
$e^{\prime}=\alpha e+a_{1}+{\cdots}+a_{r}$ with $\alpha\in R^{\prime}$ and
$a_{i}\in A_{iR^{\prime}}$. Then $e=e^{\prime}e=(\alpha e+a_{1}+{\cdots}%
+a_{r})e=\alpha e+\alpha_{1}a_{1}+{\cdots}+\alpha_{r}a_{r}$. As $\alpha_{i}\in
R^{\times}\subset R^{\prime\times}$ and $\alpha_{i}a_{i}\in A_{iR^{\prime}}$,
this implies that $a_{i}=0$ for all $i$. Therefore $e^{\prime}=\alpha e$ and
$e=\alpha e$.

(b) We have $\sigma(e)=e$ because both are left identity elements in
$A_{R^{\prime}}$. Moreover, $\sigma(A_{iR^{\prime}})$ is the submodule of
$A_{R^{\prime}}$ on which $r_{\sigma(e)}$ acts as multiplication by
$\alpha_{i}$. As $r_{\sigma(e)}=r_{e}$, we deduce that $\sigma(A_{iR^{\prime}%
})=A_{iR^{\prime}}$.
\end{proof}

\subsection{The algebra $D(L,U,S,\gamma)$}

\begin{proposition}
\label{P3}Let $V$ be a finite flat $R$-module of the form $V=L\oplus U$ with
$L$ free of rank $2$. Let $S$ be a finite flat $R$-submodule of $V^{\otimes
r}$, some $r>1$, such that $V^{\otimes r}/S$ is flat. Extend the action of
$\GL_{U}$ on $U$ to $V$ by letting it act trivially on $L$. If there exist
$\gamma_{1},\ldots,\gamma_{6}\in R$ such that the elements $0,1,\gamma
_{1},\ldots,\gamma_{6}$ are distinct modulo every maximal ideal of $R$, then
there exists a finite flat algebra $D$ over $R$ such that
\[
(\GL_{U})_{S}\simeq\underline{\Aut}(D).
\]

\end{proposition}

\begin{proof}
We define the underlying $R$-module of $D=D(L,U,S,\gamma)$ to be%
\begin{align*}
\vect(D)  &  =Re\oplus Rb\oplus Rc\oplus Rd\oplus\vect(A(V,S))\\
&  =Re\oplus Rb\oplus Rc\oplus Rd\oplus L\oplus U\oplus\bigg(\bigoplus
_{i=2}^{r-1}V^{\otimes i}\bigg)\oplus(V^{\otimes r}/S).
\end{align*}
Let $\{\ell_{1},\ell_{2}\}$ be a basis for $L$. The multiplication map on $D$
is determined by the following rules:

\begin{enumerate}
\item $e$ is a left identity element for $D;$

\item each submodule in the top row of the following table is an eigenspace
for $r_{e}$ with eigenvalue the element in the row below it,%
\[%
\begin{array}
[c]{c|c|c|c|c|c|c}%
\langle e\rangle & \langle b\rangle & \langle c\rangle & \langle d\rangle &
L & U & \left(  \bigoplus_{i=2}^{r-1}V^{\otimes i}\right)  \oplus(V^{\otimes
r}/S)\\
1 & \gamma_{1} & \gamma_{2} & \gamma_{3} & \gamma_{4} & \gamma_{5} &
\gamma_{6}%
\end{array}
;
\]

\item the multiplication table for $b,c,d$ is%

\[
\renewcommand{\arraystretch}{1.3}\arraycolsep=6pt%
\begin{array}
[c]{c|rcl}
& b & c & d\\\hline
b & 0 & c+\frac{\gamma_{2}-\gamma_{1}}{\gamma_{2}-\gamma_{3}}b & 0\\
c & -c & b & e\\
d & \ell_{1} & d & \ell_{2}.
\end{array}
\]

\item $\langle b,c,d\rangle\cdot A(V,S)=0=A(V,S)\cdot\langle b,c,d\rangle$;

\item $A(V,S)$ is a subalgebra of $D$.
\end{enumerate}

The action of $(\GL_{U})_{S}$ on $\mathrm{T}(V)_{+}$ leaves the ideal $I(S)$
stable, and so it passes to the quotient $A(V,S)$ (see the proof of
Proposition \ref{S21}). We extend this action on $\vect(A(V,S))$ to an action
on $\vect(D)$ by letting $(\GL_{U})_{S}$ act trivially on $\langle
e,b,c,d\rangle$. In this way, we get a homomorphism%
\begin{equation}
(\GL_{U})_{S}\rightarrow\underline{\Aut}(D). \label{e1}%
\end{equation}
It remains to show that this is an isomorphism.

Let $R^{\prime}$ be a commutative $R$-algebra. We have to show that the map%
\[
(\GL_{U})_{S}(R^{\prime})\rightarrow\Aut(D_{R^{\prime}})
\]
is an isomorphism. It is clearly injective. On the other hand, let $\sigma$ be
an automorphism of the algebra $D_{R^{\prime}}$ over $R^{\prime}$. According
to Lemma \ref{L1}, $\sigma(e)=e$ and $\sigma$ stabilizes each of the
$R^{\prime}$-submodules $R^{\prime}b$, $R^{\prime}c$, $R^{\prime}d$,
$L_{R^{\prime}}$, $U_{R^{\prime}}$, and $\left(  \left(  \bigoplus_{i=2}%
^{r-1}V^{\otimes i}\right)  \oplus(V^{\otimes r}/S)\right)  _{R^{\prime}}$ of
$\vect(D)_{R^{\prime}}$. Let $\sigma(b)=\gamma_{b}b$, $\sigma(c)=\gamma_{c}c$,
and $\sigma(d)=\gamma_{d}d$, where the $\gamma$ lie in $R^{\prime}$. Now%
\begin{align*}
c\cdot d=e  &  \implies\gamma_{c}\gamma_{d}=1\\
d\cdot c=d  &  \implies\gamma_{c}\gamma_{d}=\gamma_{d}\\
c\cdot b=-c  &  \implies\gamma_{c}\gamma_{b}=\gamma_{c}.
\end{align*}
From the first equation, we see that $\gamma_{c}$ and $\gamma_{d}$ are units
in $R^{\prime}$, and so the remaining equations show that $\gamma_{c}%
=1=\gamma_{b}$. Therefore $\gamma_{d}=1$ also, and so $\sigma$ acts as the
identity map on $\langle e,b,c,d\rangle_{R^{\prime}}$. Next%
\begin{align*}
d\cdot b  &  =\ell_{1}\implies\sigma(\ell_{1})=\ell_{1}\\
d\cdot d  &  =\ell_{2}\implies\sigma(\ell_{2})=\ell_{2},
\end{align*}
and so $\sigma$ acts as the identity on $L_{R^{\prime}}$. Finally, $\sigma$
acts on $\vect(A(V,S))_{R^{\prime}}$ as an automorphism of $A(V,S)_{R^{\prime
}}$. As it maps $V_{R^{\prime}}$ into $V_{R^{\prime}}$, Proposition \ref{P3}
shows that $\sigma$ arises from an element of $(\GL_{U})_{S}(R^{\prime})$.
\end{proof}

Note that $D$ is not associative: we need not have $xe\cdot y=x\cdot ey$.

\section{Algebraic groups as stabilizers}

In this section, we explain how to realize algebraic groups as the stabilizers
of submodules or of families of tensors, and we prove Theorems \ref{T3} and
\ref{T2}.

\subsection{Preliminaries}

In this subsection, we extend some standard results from base fields to more
general rings.

\begin{plain}
\label{a4}An $R$-module $V$ is finite flat if it satisfies the following
equivalent conditions (see the author's notes on commutative algebra, 12.6):

\begin{itemize}
\item $V$ is finitely generated and projective;

\item $V$ is finitely presented and flat;

\item $V$ is locally free of finite rank.
\end{itemize}

Assume that $R$ is noetherian, and let $W$ be an $R$-submodule of an
$R$-module $V$. If $V$ is finitely generated and $V/W$ is flat, then $V/W$ is
projective, and so $W$ is a direct summand of $V$, i.e., $V=W\oplus W^{\prime
}$ for some $R$-submodule $W^{\prime}$ of $V$. Conversely, if $V$ is finite
flat and $W$ is a direct summand of $V$, then $V/W$ is isomorphic to a direct
summand of $V$, and hence is (finite) flat (ibid., 11.3).
\end{plain}

\begin{plain}
\label{a29}Let $f_{1},\ldots,f_{m}\in R$ be such that $f_{1}+\cdots+f_{m}=1$.
For any $R$-module $V$,%
\[
\begin{tikzpicture}[scale=2.2,text height=1.5ex, text depth=0.25ex]
\node (a) at (0,0) {$0$};
\node (b) at (0.6,0) {$V$};
\node (c) at (1.4,0) {$\prod\nolimits_i V_{f_i}$};
\node (d) at (2.5,0) {$\prod\nolimits_{i,j} V_{f_if_j}$};
\path[->,font=\scriptsize,>=angle 90]
(a) edge (b)
(b) edge (c)
([yshift= 1pt]c.east) edge  ([yshift= 1pt]d.west)
([yshift= -1pt]c.east) edge ([yshift= -1pt]d.west);
\end{tikzpicture}
\]
is exact (ibid., 11.22). When $V$ is finite flat, the $f_{i}$ may be chosen so
that $V_{f_{i}}$ is free as an $R_{f_{i}}$-module. This often allows us in
proofs to suppose that $V$ is free.
\end{plain}

\begin{plain}
\label{a28}Assume that $R$ is an integral domain, and let $V$ and $W$ be
finite flat $R$-modules. If $v$ and $w$ are nonzero elements of $V$ and $W$,
then $v\otimes w$ is a nonzero element of $V\otimes W$. This becomes obvious
once we tensor with the field of fractions of $R$. Note that the hypothesis on
$R$ is necessary: if $R$ contains nonzero elements $a,b$ such that $ab=0$,
then $a$ and $b$ are nonzero elements of the $R$-module $R$, but $a\otimes
b=0$ in $R\otimes R\simeq R$.
\end{plain}

\begin{plain}
\label{a6}Assume that $R$ is noetherian. Let $V$ be an $R$-module, and let
$\mathrm{T}V=\bigoplus_{n}V^{\otimes n}$ be its tensor algebra. The exterior
algebra $\bigwedge V$ of $V$ is $\mathrm{T}V/I$, where $I$ is the two-sided
ideal generated by the elements $x\otimes x$, $x\in V$. The antisymmetrization
map is%
\[
a_{n}\colon V^{\otimes n}\rightarrow V^{\otimes n},\quad a_{n}(t)=\sum
_{\sigma\in S_{n}}\sign(\sigma)\sigma(t).
\]
If $V$ is finite flat, then the kernel of $a_{n}$ is $I_{n}$, and so $a_{n}$
defines an isomorphism
\[
\tstyle\bigwedge\nolimits^{n}V\rightarrow A_{n}^{\prime\prime}(V)\subset
V^{\otimes n},\quad A_{n}^{\prime\prime}(V)\overset{\df}{=}\im(a_{n});
\]
moreover, $A_{n}^{\prime\prime}(V)$ is locally a direct summand of $T^{n}V$,
and so it is finite flat. See Bourbaki, Algebra, III, \S 7, no.\thinspace4,
Remark, and Exercise 8.
\end{plain}

\begin{plain}
\label{a21}Let $V$ be a finite flat $R$-module. Then $\GL_{V}$ is a flat
algebraic group over $R$, locally isomorphic to $\GL_{n}$, $n=\rank V$.
\end{plain}

\begin{plain}
\label{a23}Let $G$ be an algebraic group over $R$ and $V$ an $R$-module. By an
action of $G$ on $V$, we mean an action of $G(R^{\prime})$ on $V(R^{\prime})$
functorial in the $R$-algebra $R^{\prime}$. When $V$ is finite flat, to give
an action of $G$ on $V$ is the same as giving a homomorphism of algebraic
groups $G\rightarrow\GL_{V}$.
\end{plain}

\begin{plain}
\label{a32}Let $V$ be a finite flat $R$-module. An action $r\colon
G\rightarrow\GL_{V}$ of $G$ on $V$ maps the universal element in
$G(\mathcal{O}{}(G))$ to an $\mathcal{O}{}(G)$-linear endomorphism of
$V\otimes\mathcal{O}{}(G)$, which is determined by its restriction to $V$,
\[
\rho\colon V\rightarrow V\otimes\mathcal{O}{}(G).
\]
The map $\rho$ is a co-action of the Hopf algebra $\mathcal{O}{}(G)$ on $V$,
i.e.,
\begin{equation}
\left\{
\begin{array}
[c]{r@{\,\,}c@{\,\,}l}%
(\id_{V}\otimes\Delta)\circ\rho & = & (\rho\otimes\id_{\mathcal{O}{}(G)}%
)\circ\rho\\
(\id_{V}\otimes\epsilon)\circ\rho & = & \id_{V}.
\end{array}
\right.  \label{e0}%
\end{equation}
In this way, we get a one-to-one correspondence $r\leftrightarrow\rho$ between
the actions of $G$ on $V$ and the co-actions of $\mathcal{O}{}(G)$ on $V$
(Milne 2017, 4.1).
\end{plain}

\begin{lemma}
\label{a30}Let $G$ be an algebraic group over $R$. Let $W$ be a finite flat
$R$-module on which $G$ acts, and let $\rho\colon W\rightarrow W_{0}%
\otimes\mathcal{O}{}(G)$ be the corresponding co-action map. Then $\rho$ is
$G$-equivariant, and realizes $W$ as a direct summand of $W_{0}\otimes
\mathcal{O}{}(G)$.
\end{lemma}

\begin{proof}
The first equality in (\ref{e0}) says that $\rho$ is a homomorphism of
$\mathcal{O}{}(G)$-comodules (and hence a homomorphism of $G$-modules). The
second equality in (\ref{e0}) says that the composite of $\rho$ with
$\id_{V_{0}}\otimes\epsilon$ is the identity map.
\end{proof}

\begin{plain}
\label{a22}Let $G$ be an algebraic group over $R$ and $V$ an $R$-module on
which $G$ acts. When $i\colon S\rightarrow V$ is an $R$-\textit{submodule} of
$V$, we define $G_{S}$ (stabilizer of $S$ in $G$) to be the functor%
\[
R^{\prime}\rightsquigarrow\{\alpha\in\Aut_{R^{\prime}}(V_{R^{\prime}}%
)\mid\alpha(i_{R^{\prime}}(S_{R^{\prime}}))=i_{R^{\prime}}(S_{R^{\prime}})\}.
\]
When $S$ is a \textit{subset} of $V$, we define $G_{S}$ to be the functor%
\[
R^{\prime}\rightsquigarrow\{\alpha\in\Aut_{R^{\prime}}(V_{R^{\prime}}%
)\mid\alpha(s)=s\text{ for all }s\in S\}.
\]
If $S$ is contained in an $R$-submodule $V^{\prime}$ of $V$, stable under $G$,
and $V/V^{\prime}$ is flat, then the group functor $G_{S}$ is the same for
$S\subset V^{\prime}$ as for $S\subset V$ (because the map $V_{R^{\prime}%
}^{\prime}\rightarrow V_{R^{\prime}}$ is injective for all $R$-algebras
$R^{\prime}$).
\end{plain}

\begin{plain}
\label{a5}When $R$ is noetherian, every comodule over a flat $R$-coalgebra is
a filtered union of finitely generated subcomodules (Serre 1993, 1.4). In
particular, every $G$-module, where $G$ is a flat algebraic group over $R$, is
a filtered union of finite generated $G$-submodules.
\end{plain}

\begin{lemma}
\label{a7}Let $R$ be an integral domain and $G$ an algebraic group over $R$.

\begin{enumerate}
\item Let $V_{1}$ and $V_{2}$ be finite flat $R$-modules on which $G$ acts,
and let $S_{1}\subset V_{1}$ and $S_{2}\subset V_{2}$ be nonzero submodules
such that $V_{1}/S_{1}$ and $V_{2}/S_{2}$ are flat. Then the stabilizer of
$S_{1}\otimes S_{2}\subset V_{1}\otimes V_{2}$ in $\GL_{V_{1}}\times
\GL_{V_{2}}$ is equal to the stabilizer of $S_{1}\oplus S_{2}\subset
V_{1}\oplus V_{2}$ in $\GL_{V_{1}}\times\GL_{V_{2}}$.

\item Let $V_{1}=V=V_{2}$ in (a). Then the stabilizer of $S_{1}\otimes
S_{2}\subset V\otimes V$ in $\GL_{V}$ is equal to the intersection of the
stabilizers of $S_{1}\subset V$ and $S_{2}\subset V$ in $\GL_{V}$.

\item Let $V$ be a finite flat $R$-module on which $G$ acts, and let $L$ be a
line (i.e., one-dimensional subspace) in $V$ such that $V/L$ is flat. For
every $r>0$, the stabilizer of $L\subset V$ in $\GL_{V}$ is equal to the
stabilizer of $L^{\otimes r}\subset V^{\otimes r}$ in $\GL_{V}$.
\end{enumerate}
\end{lemma}

\begin{proof}
(a) As $V_{1}/S_{1}$ and $V_{2}/S_{2}$ are flat and finitely generated and $R$
is integral domain, they are finitely presented. The hypotheses imply that
$V_{1}=S_{1}\oplus W_{1}$ and $V_{2}=S_{2}\oplus W_{2}$ for some finite flat
$R$-submodules $W_{1}$ and $W_{2}$ of $V_{1}$ and $V_{2}$. Then%
\[
V_{1}\otimes V_{2}=\left(  S_{1}\otimes S_{2}\right)  \oplus\left(
S_{1}\otimes W_{2}\right)  \oplus\left(  W_{1}\otimes S_{2}\right)
\oplus\left(  W_{1}\otimes W_{2}\right)  \text{.}%
\]
Let $R^{\prime}$ be an $R$-algebra and $\alpha_{1}$ and $\alpha_{2}$
automorphisms of $V_{1R^{\prime}}$ and $V_{2R^{\prime}}$. We have to show that%
\[
(\alpha_{1}\otimes\alpha_{2})(S_{1}\otimes S_{2})\subset S_{1}\otimes
S_{2}\implies\alpha_{1}(S_{1})\subset S_{1}\text{ and }\alpha_{2}%
(S_{2})\subset S_{2},
\]
the reverse implication being obvious. Let $s_{1}$ and $s_{2}$ be nonzero
elements of $S_{1}$ and $S_{2}$, and let $\alpha_{1}(s_{1})=s_{1}^{\prime
}+w_{1}$ and $\alpha_{2}(s_{2})=s_{2}^{\prime}+w_{2}$. Then%
\[
S_{1}\otimes S_{2}\ni(\alpha_{1}\otimes\alpha_{2})(s_{1}\otimes s_{2}%
)=s_{1}^{\prime}\otimes s_{2}^{\prime}+s_{1}^{\prime}\otimes w_{2}%
+w_{1}\otimes s_{2}^{\prime}+w_{1}\otimes w_{2}.
\]
If $w_{1}\neq0$, then $s_{2}^{\prime}=0=w_{2}$ (see \ref{a28}), contradicting
$s_{2}\neq0$. Hence $w_{1}=0$, and similarly, $w_{2}=0$.

Statement (b) follows from (a), and (c) follows from (b).

\begin{lemma}
\label{a8}Assume that $R$ is noetherian. Let $V$ be a finite flat $R$-module
and $S$ an $R$-submodule such that $V/S$ is flat. If $S$ is locally free of
rank $d$, then the stabilizer of $S\subset V$ in $\GL_{V}$ is equal to the
stabilizer of $\bigwedge^{d}S\subset\bigwedge^{d}V$ in $\GL_{V}$.\footnote{This is a standard fact, implicit in the proof of
the projectivity of Grassmanians.}
\end{lemma}

\begin{proof}
If $S=V$, this is obvious, and so we assume that $S\neq V$. Because $V/S$ is
flat, $V=S\oplus W$ for some $R$-submodule $W$ of $V$ (here we use that $R$ is
noetherian). As $V$ is finite flat, so also is $W$.

Let $L=\bigwedge^{d}S$. Let $R^{\prime}$ be an $R$-algebra and $\alpha$ an
automorphism of $V_{R^{\prime}}$. We shall show that%
\[
\alpha L_{R^{\prime}}=L_{R^{\prime}}\implies\alpha S_{R^{\prime}}%
=S_{R^{\prime}},
\]
the reverse implication being obvious. We may suppose that $S$ and $W$ are
free (\ref{a29}).

Let $(e_{j})_{1\leq i\leq d}$ be a basis for $S$, and extend it to a basis
$(e_{i})_{1\leq i\leq n}$ of $V$. Let $s=e_{1}\wedge\cdots\wedge e_{d}$. Then%
\[
\tstyle S_{R}=\{v\in V_{R}\mid s\wedge v=0\text{ (in }\bigwedge\nolimits^{d+1}%
V_{R}\text{)}\}\text{.}%
\]
To see this, let $v=\sum_{i=1}^{n}a_{i}e_{i}$, $a_{i}\in R$, be an element of
$V_{R}$. Then%
\[
s\wedge v=\sum_{d+1\leq i\leq n}a_{i}e_{1}\wedge\cdots\wedge e_{d}\wedge
e_{i}\text{.}%
\]
As the elements $e_{1}\wedge\cdots\wedge e_{d}\wedge e_{i}$, $d+1\leq i\leq
n$, are part of a basis for $\bigwedge\nolimits^{d+1}V$, we see that
\[
s\wedge v=0\iff a_{i}=0\text{ for all }d+1\leq i\leq n\iff v\in S\text{.}%
\]

Let $\alpha\in\GL(V_{R})$. If $(\bigwedge\nolimits^{d}\alpha)(L_{R})=L_{R}$,
then $(\bigwedge\nolimits^{d}\alpha)s=cs$ for some $c\in R^{\times}$. If $v\in
S_{R}$, then $s\wedge v=0$, and so
\[
\tstyle0=(\bigwedge\nolimits^{d+1}\alpha)(s\wedge v)=(\bigwedge\nolimits^{d}%
\alpha)s\wedge\alpha v=c\left(  s\wedge(\alpha v)\right)  ,
\]
which implies that $\alpha v\in S_{R}$.
\end{proof}

Recall that there is a natural left action of $G$ on $\mathcal{O}{}(G)$ (the
regular representation), namely,
\[
(gf)(x)=f(xg),\quad f\in\mathcal{O}{}(G),\quad g\in G,\quad x\in G.
\]

\begin{lemma}
\label{a19}Let $G$ be an algebraic group over $R$ and $H$ a closed algebraic
subgroup of $G$. Let $I\subset\mathcal{O}{}(G)$ be the ideal of $H$. Then $H$
is the stabilizer of $I$ in $\mathcal{O}(G)$, i.e., for all $R$-algebras
$R^{\prime}$,%
\[
H(R^{\prime})=\{g\in G(R^{\prime})\mid gI_{R^{\prime}}\subset I_{R^{\prime}%
}\}.
\]

\end{lemma}

\begin{proof}
Let $h\in H(R^{\prime})$ some $R^{\prime}$, and let $f\in I_{R^{\prime}}$.
Then, for all $R^{\prime}$-algebras $R^{\prime\prime}$ and $x\in
H(R^{\prime\prime})$,
\[
(hf)(x)\overset{\df}{=}f(xh)=0
\]
because $xh\in H(R^{\prime\prime})$. Hence $hf\in I_{R^{\prime}}$.

Let $g\in G(R^{\prime})$ be such that $gI_{R^{\prime}}\subset I_{R^{\prime}}$,
and let $f\in I$. Then%
\[
f(g)=f(e\cdot g)=(gf)(e)=0\text{,}%
\]
because $gf\in I_{R^{\prime}}$. Hence $g\in H(R^{\prime})$.
\end{proof}

\begin{plain}
\label{a24}Let $V$ be a finite flat $R$-module. We let $\GL_{V}$ act on the
(finite flat) $R$-module $\End(V)$ by setting%
\[
g\alpha=g\circ\alpha,\quad g\in\GL(V_{R^{\prime}}),\quad\alpha\in
\End(V_{R^{\prime}}),\quad\text{some }R^{\prime}.
\]
Then the canonical isomorphism $\End(V)\simeq V_{0}^{\vee}\otimes V$ of
$R$-modules is $\GL_{V}$-equivariant. To check this, let $f\otimes v\in
V_{0}^{\vee}\otimes V$, and regard it as the element of $\End(V)$ such that%
\[
(f\otimes v)(x)=f(x)v,\quad x\in V.
\]
For $g\in\GL(V)$,
\[
(g(f\otimes v))(x)=g((f\otimes v)(x))=g(f(x)v)=f(x)gv=(f\otimes gv)(x),
\]
and so $g(f\otimes v)=(f\otimes gv)$ as claimed.
\end{plain}

\begin{plain}
\label{a20}Let $V$ be a finite flat $R$-module and $G$ a closed algebraic
subgroup of $\GL_{V}$. Then $\GL_{V}$ is a schematically dense open subscheme
of $\End_{V}$ (multiplicative monoid scheme). Correspondingly
\[
\Sym(\End(V))=\mathcal{O}{}(\End_{V})\subset\mathcal{O}{}(\GL_{V})\text{.}%
\]
For example, if $V$ is free, then the choice of a basis for $V$ identifies the
inclusion with%
\[
R[X_{ij}]\subset R[X_{ij}][1/\det(X_{ij})].
\]
The inclusion $\Sym(\End(V))\hookrightarrow\mathcal{O}{}(\GL_{V})$ is
$\GL_{V}$-equivariant for the actions considered in \ref{a24}. Let $I$ be the
ideal of $G$ in $\mathcal{O}{}(\GL_{V})$, and let $I^{\prime}=I\cap
\Sym(\End(V))$. Then $I^{\prime}$ generates the ideal $I$, and so $G$ is the
stabilizer of $I^{\prime}\subset\Sym(\End(V))$ in $\GL_{V}$ (Lemma \ref{a19}).
\end{plain}
\end{proof}

\begin{plain}
\label{a11}Recall that an algebraic group over a field is said to be linearly
reductive if every finite-dimensional representation is semisimple. In
characteristic zero, $G$ is linearly reductive if and only if $G^{\circ}$ is
reductive. In characteristic $p\neq0$, $G$ is linearly reductive if and only
if $G^{\circ}$ is of multiplicative type and $p$ does not divide the index
$(G\colon G^{\circ})$ (Nagata's theorem). See Milne 2017, 12.56.
\end{plain}

\subsection{Are algebraic groups linear?}

Let $G$ be a flat algebraic group over a ring $R$. Does there exist an
embedding of $G$ into $\GL_{n}$ for some $n$? Apparently the answer is not
known even for $R$ the ring of dual numbers over a field. However, there is
the following result.\footnote{This should
be considered folklore. See an earlier footnote.}

\begin{proposition}
\label{a31}Let $G$ be a flat algebraic group over a Dedekind domain $R$. There
exists a finite flat $R$-submodule $V$ of $\mathcal{O}(G)$, stable under $G$,
such that the homomorphism $G\rightarrow\GL_{V}$ is a closed immersion.
\end{proposition}

\begin{proof}
There exists a finitely generated $R$-submodule $V$ of $\mathcal{\mathcal{O}%
{}}(\GL_{V})$, stable under $G$, containing a set of generators for
$\mathcal{O}{}(G)$ (see \ref{a5}). Now $G$ flat over $R\implies\mathcal{O}%
{}(G)$ is torsion-free (as an $R$-module)$\implies V$ is torsion-free$\implies
V$ is flat (because $R$ is a Dedekind domain). It remains to show that the
homomorphism $\alpha\colon\mathcal{O}(\GL_{V})\rightarrow\mathcal{O}(G)$
defined by the action of $G$ on $V$ is surjective.

We may suppose that $V$ is free (see \ref{a29}). Let $\Delta\colon
\mathcal{O}{}(G)\rightarrow\mathcal{O}{}(G)\otimes\mathcal{O}{}(G)$ be the
comultiplication map and $\epsilon\colon\mathcal{O}{}(G)\rightarrow R$ the
co-identity map. Let $(e_{i})_{1\leq i\leq n}$ be a basis for $V$, and write
$\Delta(e_{j})=\sum_{i}e_{i}\otimes a_{ij}$, $a_{ij}\in\mathcal{O}{}(G)$. The
image of $\alpha$ contains the $a_{ij}$ (the choice of the basis $(e_{i})$,
determines an isomorphism $\mathcal{O}{}(\GL_{V})\simeq R[T_{ij}]$, and
$\alpha$ maps $T_{ij}$ to $a_{ij}$; see Milne 2017, 4.1). As $\epsilon
\colon\mathcal{O}{}(G)\rightarrow R$ is the co-identity,
\[
e_{j}=(\epsilon\otimes\id_{A})\Delta(e_{j})=\sum\nolimits_{i}\epsilon
(e_{i})a_{ij},
\]
and so the image of $\alpha$ contains $V$, which we chose to generate
$\mathcal{O}{}(G)$.
\end{proof}

Thus, if $R$ is a Dedekind domain and $G$ is flat, then there is an embedding
of $G$ into $\GL_{V}$ for some finite flat $R$-module $V$. Such a $V$ is a
direct summand $F=V\oplus W$ of a free $R$-module $F$ of finite rank. Extend
the action of $G$ on $V$ to $F$ by letting it act trivially on $W$ and choose
a basis for $F$. Now $G$ is a closed algebraic subgroup of $\GL_{n}$,
$n=\rank
F$.

\subsection{Expressing all representations in terms of one faithful
representation}

Let $G$ be an algebraic group over a field $k$, and let $(V,r)$ be a faithful
representation of $G$. Then $V$ generates the tannakian category of
finite-dimensional representations of $G$. This means that every
finite-dimensional representation of $G$ can be constructed from $V$ by
forming tensor products, direct sums, duals, and subquotients (Milne 2017, 4.14).

In this section, we present variants of this statement. For a finite flat
$R$-module of rank $r$, we let $\det=\bigwedge^{r}V$ and $\det^{-1}=\det
^{\vee}$. For an $R$-module $V$, we let $T^{m,n}(V)=V^{\otimes m}%
\otimes(V^{\vee})^{\otimes n}$.

\begin{proposition}
\label{a18}Let $G$ be a closed algebraic subgroup of $\GL_{V}$, where $V$ is a
free $R$-module of finite rank. Let $W$ be a $G$-module that is free of finite
rank as an $R$-module. For some $s$, $W\cdot\det^{s}$ is isomorphic to a
submodule of a quotient of a direct sum of tensor powers of $V$.
\end{proposition}

\begin{proof}
The choice of a basis for $W_{0}$ realizes $W$ as a $G$-submodule of
$\mathcal{O}(G)^{m}$, $m=\rank W$ (see Lemma \ref{a30}). The embedding
$G\hookrightarrow\GL_{V}$ corresponds to a surjective homomorphism
$\mathcal{O}{}(\GL_{V})\rightarrow\mathcal{O}{}(G)$. Recall that
$\mathcal{O}{}(\GL_{V})=\Sym(\End(V))[1/\det]$ and that $\End(V)\simeq
V_{0}^{\vee}\otimes V$ as a $G$-module (\ref{a24}). The choice of a basis for
$V_{0}$ determines a $G$-isomorphism $\End(V)\simeq nV$, $n=\rank V$. We have
$G$-homomorphisms%
\[
\mathrm{T}(nV)^{m}\twoheadrightarrow\Sym(nV)^{m}\subset\mathcal{O}{}%
(\GL_{V})^{m}\twoheadrightarrow\mathcal{O}{}(G)^{m}.
\]
For some $s\geq0$, $W\cdot\det^{s}$ is contained in the image of
$\mathrm{T}(nV)^{m}$ in $\mathcal{O}{}(G)^{m}$ . Hence $W\cdot\det^{s}$ is
contained in a quotient of $\mathrm{T}^{\leq h}(nV)^{m}$ for some $h$, and
$\mathrm{T}^{\leq h}(nV)^{m}$ is a sum of tensor powers of $V$.
\end{proof}

\setcounter{corollary}{0}

\begin{corollary}
\label{a33}Let $G$, $V$, and $W$ be as in the proposition. Then $W$ is
isomorphic to a submodule of a quotient of a direct sum of modules
$T^{m,n}(V)$.
\end{corollary}

\begin{proof}
Let $n=\rank V$. As $\det$ is a direct summand of $V^{\otimes n}$ (see
\ref{a6}), its dual $\det^{-1}$ is a direct summand of $\left(  V^{\vee
}\right)  ^{\otimes n}$. In the proof of Proposition \ref{a18}, we constructed
a diagram%
\[
W\otimes\det\nolimits^{s}\hookrightarrow Q\twoheadleftarrow\mathrm{T}^{\leq
h}(nV)^{m}.
\]
On tensoring this with $(V^{\vee})^{\otimes ns}$, we get a diagram%
\[
W\subset W\otimes\det\nolimits^{s}\otimes(V^{\vee})^{\otimes ns}%
\hookrightarrow Q\otimes(V^{\vee})^{\otimes ns}\twoheadleftarrow
\mathrm{T}^{\leq h}(nV)^{m}\otimes(V^{\vee})^{\otimes ns},
\]
as required.
\end{proof}

\begin{remark}
\label{a35}If $R$ is a field and $G$ is linearly reductive, then
\textquotedblleft of a quotient\textquotedblright\ can be omitted from the
statements of Proposition \ref{a18} and Corollary \ref{a33}.
\end{remark}

When $G\subset\SL_{V}$, the proof of Proposition \ref{a18} simplifies.

\begin{proposition}
\label{a1}Let $G$ be a closed algebraic subgroup of $\SL_{V}$, where $V$ is a
free $R$-module of finite rank. Let $W$ be a $G$-module that is free of finite
rank as an $R$-module. Then $W$ is isomorphic to a submodule of a quotient of
a direct sum of tensor powers of $V$.
\end{proposition}

\begin{proof}
As before, $W\subset\mathcal{O}{}(G)^{m}$, $m=\rank W$. In this case, we get
$G$-homomorphisms%
\[
T(nV)^{m}\twoheadrightarrow\Sym(nV)^{m}\twoheadrightarrow\mathcal{O}%
(\SL_{V})^{m}\twoheadrightarrow\mathcal{O}{}(G)^{m}.
\]
For some $h$, $W$ is contained in a quotient of $T^{\leq h}(nV)^{m}$.
\end{proof}

When $V$ is a finite-dimensional vector space over a field $k$ of
characteristic zero, Proposition \ref{a1} shows that every finite-dimensional
$\SL_{V}$-module $W$ is isomorphic to a submodule of $T(nV)^{m}$, where
$n=\dim V$ and $m=\dim W$. In fact, a stronger result holds.

\begin{proposition}
[Gordeev--Popov]\label{a16}Let $V$ be a finite-dimensional vector space over a
field $k$. Every finite-dimensional $\SL_{V}$-module is isomorphic to a
submodule of $T(V)_{+}$.
\end{proposition}

\begin{proof}
See Gordeev and Popov 2003, Proposition 11.
\end{proof}

\subsection{Algebraic groups as stabilizers}

\begin{proposition}
\label{a39}Let $R$ be a noetherian ring. Let $G$ be a closed algebraic
subgroup of $\GL_{V}$, where $V$ is a finite flat $R$-module. For some
$h\geq0$, there exists an $R$-submodule $S\subset\mathrm{T}^{\leq h}%
(V_{0}^{\vee}\otimes V)$ such that $G$ is the stabilizer of $S$ in $\GL_{V}$.
\end{proposition}

\begin{proof}
Let $I$ be the kernel of the homomorphism of $R$-algebras%
\[
\Sym(V_{0}^{\vee}\otimes V)\hookrightarrow\mathcal{O}{}(\GL_{V}%
)\twoheadrightarrow\mathcal{O}{}(G).
\]
Then $G$ is the stabilizer of $I$ in $\GL_{V}$ (see \ref{a20}). For some
$h\geq0$, $\Sym^{\leq h}(V_{0}^{\vee}\otimes V)$ contains a set of generators
for the ideal $I$ (here we use that $R$ is noetherian), and $G$ is the
stabilizer of
\[
I\cap\Sym^{\leq h}(V_{0}^{\vee}\otimes V)\subset\Sym^{\leq h}(V_{0}^{\vee
}\otimes V)
\]
in $\GL_{V}$. Now $G$ is the stabilizer in $\GL_{V}$ of the preimage $S$ of
$I\cap\Sym^{\leq h}(V_{0}^{\vee}\otimes V)$ under the quotient map%
\[
\mathrm{T}^{\leq h}(V_{0}^{\vee}\otimes V)\twoheadrightarrow\Sym^{\leq
h}(V_{0}^{\vee}\otimes V)\text{.}%
\]

\end{proof}

\begin{remark}
\label{a43}If $R$ is a Dedekind domain and $G$ is flat, then the $R$-submodule
$S$ constructed in the proof of the proposition has the property that
$\mathrm{T}^{\leq h}(V_{0}^{\vee}\otimes V)/S$ is flat. To see this, note that
the hypotheses imply that $\Sym(V_{0}^{\vee}\otimes V)/I$ is torsion-free, and
so $I$ is saturated as an $R$-submodule of $\Sym(V_{0}^{\vee}\otimes V)$. It
follows that $I\cap\Sym^{\leq h}(V_{0}^{\vee}\otimes V)$ is saturated, and so%
\[
\mathrm{T}^{\leq h}(V_{0}^{\vee}\otimes V)/S\simeq\Sym^{\leq h}(V_{0}^{\vee
}\otimes V)/I\cap\Sym^{\leq h}(V_{0}^{\vee}\otimes V)
\]
is flat.
\end{remark}

The next statement improves results of Deligne (1982, 3.1) and Kisin (2010,
1.3.1). It has applications to Shimura varieties in mixed characteristic.

\begin{proposition}
\label{a40}Let $R$ be a Dedekind domain. Let $G$ be a closed algebraic
subgroup of $\GL_{V}$, where $V$ is a finite flat $R$-module. If the generic
fibre of $G$ is linearly reductive, then, locally on $\Spec R$, $G$ is the
algebraic subgroup of $\GL_{V}$ fixing a finite collection of tensors in
spaces $V^{\otimes m}\otimes(V^{\vee})^{\otimes n}$, $m,n\geq0$.
\end{proposition}

\begin{proof}
By \textquotedblleft locally on $\Spec R$\textquotedblright\ we mean that
there exist $f_{i}\in R$ such that $f_{1}+\cdots+f_{m}=1$ and the statement
holds after a base change $R\rightarrow R_{f_{i}}$. Thus, we may suppose that
$V$ is free$,$ say, of rank $n$, and replace $V_{0}^{\vee}\otimes V$ with $nV$
in Proposition \ref{a39}. Let $S\subset T^{\leq h}(nV)\overset{\df}{=}W$ be as
in that proposition. Then $W$ is free of finite rank, and so $S$ is finite
flat (here we use that $R$ is Dedekind). Let $r=\rank S$. Then $G$ is the
stabilizer of $L\overset{\df}{=}\bigwedge^{r}S\subset\bigwedge^{r}W$ in
$\GL_{V}$ (Lemma \ref{a8}). Note that $L$ is locally free of rank $1$ and that
$\bigwedge^{r}W$ is a direct summand of $\bigotimes\nolimits^{r}W$ (see Lemma
\ref{a15}), which is a direct sum of tensor powers of $V$.

As the generic fibre of $G$ is linearly reductive, the quotient map $\left(
\bigwedge^{r}W\right)  ^{\vee}\rightarrow L^{\vee}$ has a $G$-equivariant
section over the generic point $\eta$ of $\Spec R$. It follows that there
exists a $G$-stable line $L^{\ast}\subset\left(  \bigwedge^{r}W\right)
^{\vee}$ that maps isomorphically to $L^{\vee}$ over $\eta$. Now $G$ acts
trivially on $L\otimes_{R}L^{\ast}$ because this is so over $\eta$, and the
stabilizer of
\[
L\otimes_{R}L^{\ast}\subset\bigwedge\nolimits^{r}W\otimes(\bigwedge
\nolimits^{r}W^{\vee})\subset\bigotimes\nolimits^{r}W\otimes(\bigotimes
\nolimits^{r}W^{\vee})
\]
in $\GL_{V}$ is equal to $G$.

After a base change $R\rightarrow R_{f_{i}},$ the module $L\otimes_{R}L^{\ast
}$ will be free. Let $\{s\}$ be a basis for $L\otimes_{R}L^{\ast}$, and write
$s=\sum_{i\in I}s_{i}$ with each $s_{i}$ an element of a module $T^{m,n}(V)$.
Then $G=(\GL_{V})_{S}$ with $S=\{s_{i}\mid i\in I\}.$
\end{proof}

\setcounter{corollary}{0}

\begin{corollary}
\label{a41}Let $G$ be a flat algebraic group over a Dedekind domain $R$. There
exists a finite flat $R$-module $V$ and an embedding $G\hookrightarrow\GL_{V}%
$. If $R$ is a principal ideal domain and the generic fibre of $G$ over $R$ is
linearly reductive, then $G$ is the algebraic subgroup of $\GL_{V}$ fixing a
finite collection of tensors in spaces $V^{\otimes m}\otimes(V^{\vee
})^{\otimes n}$, $m,n\geq0$.
\end{corollary}

\begin{proof}
This follows from Propositions \ref{a31} and \ref{a40}.
\end{proof}

\begin{remark}
\label{a14}The condition that $G_{\eta}$ is linearly reductive can be replaced
by the following condition: the map $\Hom_{\eta}(G_{\eta},\mathbb{G}%
_{m})\rightarrow\Hom_{\eta}(\GL_{U_{\eta}},\mathbb{G}_{m})$ has finite
cokernel. The proof requires Lemma \ref{a7}(c).
\end{remark}

\subsection{Algebraic groups as automorphism groups of algebras}

The next two lemmas are adapted from Gordeev and Popov 2003.\footnote{Readers
should be careful not to confuse the tensor algebra with the symmetric
algebra, for which the proof of Lemma 7 fails.}

\begin{lemma}
\label{a3}Let $U$ (resp.~$L$) be a free $R$-module of finite rank $m$
(resp.\thinspace rank $1$). There exists an injective homomorphism of graded
$\GL_{U}$-modules
\[
\iota\colon\mathrm{T}(mU)\hookrightarrow\mathrm{T}(L\oplus U)
\]
realizing $\mathrm{T}(mU)$ as a direct summand of $\mathrm{T}(L\oplus U)$.
Here $\GL_{U}$ acts trivially on $L$.
\end{lemma}

\begin{proof}
Let $U_{i}$ be the $i$th summand of $mU$ considered as a subspace of $mU$, and
choose a basis $\{f_{ij}\mid j=1,\ldots,m\}$ of $U_{i}$. Let \{$l\}$ be a
basis for $L$, and set
\[
\iota(f_{i_{1}j_{1}}\otimes\cdots\otimes f_{i_{t}j_{t}})=l^{\otimes i_{1}%
}\otimes f_{i_{1}j_{1}}\otimes\cdots\otimes l^{\otimes i_{t}}\otimes
f_{i_{t}j_{t}}.
\]
The map $\mathrm{T}(mU)\rightarrow\mathrm{T}(L\oplus U)$, defined on a basis
of $\mathrm{T}(mU)_{+}$ by this formula and sending $1$ to $1$, has the
claimed properties.
\end{proof}

When $R$ is a field, there even exists an injective homomorphism
$\mathrm{T}^{\leq h}(mU)\hookrightarrow\mathrm{T}_{+}(U)$ (Proposition
\ref{a16}).

\begin{lemma}
\label{a15} Let $U$ be a finite flat $R$-module and $L$ a free $R$-module of
rank $1$. For all $r\geq h\geq2,$ there is an injective homomorphism of
$\GL_{U}$-modules%
\[
\mathrm{T}^{\leq h}(U)\hookrightarrow(L\oplus U)^{\otimes r}%
\]
realizing $\mathrm{T}^{\leq h}(U)$ as a direct summand of $(L\oplus
U)^{\otimes r}$.
\end{lemma}

\begin{proof}
For any $r\geq1,$%
\begin{align*}
(L\oplus U)^{\otimes r}  & \simeq\bigoplus_{i+j=r}L^{\otimes i}\otimes
U^{\otimes j}\oplus\text{other terms}\\
& \simeq\mathrm{T}^{\leq r}U\oplus\text{other terms}%
\end{align*}
(the second isomorphism depends on a choice of a basis for $L$).
\end{proof}

\begin{proposition}
\label{a38}Let $R$ be a Dedekind domain. Let $G$ be a closed algebraic
subgroup of $\GL_{U}$ flat over $R$, where $U$ is a free $R$-module of finite
rank. Let $L$ be a free $R$-module of rank $2$ with $G$ acting trivially.
There exists a finite flat $R$-module $S$ of $\left(  L\oplus U\right)
^{\otimes r}$, some $r\geq2$, such that $\left(  L\oplus U\right)  /S$ is flat
and $G=(\GL_{U})_{S}$.
\end{proposition}

\begin{proof}
Let $m=\rank U$. According to Proposition \ref{a39} and Remark \ref{a43},
$G=(\GL_{U})_{S}$ with $S$ a finite flat $R$-submodule of $\mathrm{T}^{\leq
h}(mU)$ such that $\mathrm{T}^{\leq h}(mU)/S$ flat. According to Lemmas
\ref{a3} and \ref{a15}, for all $r\geq h$, there exists an injective
homomorphism $\mathrm{T}^{\leq h}(mU)\hookrightarrow(L\oplus U)^{\otimes r}$
making $\mathrm{T}^{\leq h}(mU)$ a direct summand of $(L\oplus U)^{\otimes r}$.

On combining the last two statements, we find that $G=(\GL_{U})_{S}$, where
$S$ is a finite flat $R$-submodule of $(L\oplus U)^{\otimes r}$ such that
$(L\oplus U)^{\otimes r}/S$ is flat.
\end{proof}

\begin{theorem}
\label{a36}Let $G$ be an algebraic group flat over a Dedekind domain $R$. If
there exist $\gamma_{1},\ldots,\gamma_{6}\in R$ such that the elements
$0,1,\gamma_{1},\ldots,\gamma_{6}$ are distinct modulo every maximal ideal of
$R$, then there exists a finite flat algebra $D$ over $R$ such that
$G=\underline{\Aut}(D)$.
\end{theorem}

\begin{proof}
According to Proposition 3, there exists a finite flat $R$-module $U$ and an
embedding $G\rightarrow\GL_{U}$. Now we can apply Proposition \ref{a38} and
Proposition \ref{P3}.
\end{proof}

Theorem \ref{a36} leaves open the question: given an algebraic group $G$ over
$R$, what can be said about the algebras $A$ over $R$ such that
$G=\underline{\Aut}(A)$. When $R$ is a field, Gordeev and Popov (2003) prove a
number of results about this, for example, that $A$ can be chosen to be simple.

\vspace{5pt}\centerline{\LARGE\bfseries{References}}\vspace{10pt}
\parindent 0pt \everypar{\hangindent1.5em\hangafter1}

Borel, Armand. 1991. Linear algebraic groups. Second edition. Graduate Texts
in Mathematics, 126. Springer-Verlag, New York.

Deligne, P. 1982. Hodge cycles on abelian varieties (notes by J.S. Milne). In:
Hodge Cycles, Motives, and Shimura Varieties, Lecture Notes in Math. 900,
Springer-Verlag, pp. 9--100.

Gordeev, Nikolai L.; Popov, Vladimir L. 2003. Automorphism groups of finite
dimensional simple algebras. Ann. of Math. (2) 158, no. 3, 1041--1065.

Kisin, Mark. 2010. Integral models for Shimura varieties of abelian type. J.
Amer. Math. Soc. 23, no. 4, 967--1012.

Milne, J. S. 1994. Shimura varieties and motives. Motives (Seattle, WA, 1991),
447--523, Proc. Sympos. Pure Math., 55, Part 2, Amer. Math. Soc., Providence, RI.

Milne, J. S. 2009. Rational Tate classes. Mosc. Math. J. 9, no. 1, 111--141.

Milne, J. S. 2017. Algebraic groups. The theory of group schemes of finite
type over a field. Cambridge Studies in Advanced Mathematics, 170. Cambridge
University Press, Cambridge (corrected reprint 2022).

Perepechko, Alexander. 2009. Affine algebraic monoids as endomorphisms'
monoids of finite-dimensional algebras. Proc. Amer. Math. Soc. 137, no. 10, 3227--3233.

Saavedra Rivano, Neantro. 1972. Cat\'egories Tannakiennes. Lecture Notes in
Mathematics, Vol. 265. Springer-Verlag, Berlin-New York.

Serre, Jean-Pierre. 1993. G\`{e}bres. Enseign. Math. (2) 39, no. 1-2, 33--85.

\parindent 10pt \everypar{\hangindent0em\hangafter0}

\printendnotes

\end{document}